\documentclass[10pt,twoside]{article}
\usepackage[latin1]{inputenc}
\usepackage{amsmath}
\usepackage{dsfont}
\usepackage{graphicx}
\usepackage{yhmath}
\usepackage[table]{xcolor}
\usepackage{mathrsfs} % permet d'avoir le D cursive
\usepackage{amssymb}
\usepackage{url}
\usepackage{makecell}
\usepackage{arydshln}
\usepackage{multirow}
\usepackage{arydshln}
\usepackage{mathtools}
\usepackage{stmaryrd}  % pour intervalle entiers

\input epsf
\def\limiten{\renewcommand{\arraystretch}{0.5}
\begin{array}[t]{c}\stackrel{}{\longrightarrow} \\
{\scriptstyle n\rightarrow
\infty}\end{array}\renewcommand{\arraystretch}{1}}

\def\limiteproban{\renewcommand{\arraystretch}{0.5}
\begin{array}[t]{c}\stackrel{{\cal P}}{\longrightarrow} \\
{\scriptstyle n\rightarrow
\infty}\end{array}\renewcommand{\arraystretch}{1}}

\numberwithin{equation}{section}

\newtheorem{thm}{Theorem}[section]

\newtheorem{Def}[thm]{Definition\rm}

\newtheorem{prop}[thm]{Proposition}

\newtheorem{rmrk}[thm]{Remark}

\newcommand{\E}{\ensuremath{\mathbb{E}}}
\newcommand{\R}{\ensuremath{\mathbb{R}}}
\newcommand{\Z}{\ensuremath{\mathbb{Z}}}

\newcommand{\N}{\ensuremath{\mathbb{N}}}

\newcommand{\cov}{\ensuremath{\mathrm{Cov}}}

\newcommand{\lip}{\ensuremath{\mathrm{Lip}}}
\definecolor{grisclair}{gray}{0.9}

\DeclareMathOperator*{\argmin}{argmin}

\newcommand{\schi}{\ensuremath{\mbox{\footnotesize$\chi$}}}

\newcommand{\bchi}{\ensuremath{\mbox{\large$\chi$}}}

\renewcommand{\arraystretch}{.8}
   
\newcommand{\qed}{\hfill $\blacksquare$}

\begin{document}
%\date{ }
\title{\bf Statistical learning for $\psi$-weakly dependent processes}
 \maketitle \vspace{-1.0cm}
\begin{center}
   Mamadou Lamine DIOP \footnote{Supported by
   % the Institute for advanced studies - IAS (Université de Cergy-Pontoise, France), 
   the MME-DII center of excellence (ANR-11-LABEX-0023-01) 
   %and by the CEA-MITIC (Université Gaston Berger, Sénégal)
   } 
   and 
     William KENGNE \footnote{Developed within the ANR BREAKRISK: ANR-17-CE26-0001-01 and the  CY Initiative of Excellence (grant "Investissements d'Avenir" ANR-16-IDEX-0008), Project "EcoDep" PSI-AAP2020-0000000013} 
 \end{center}

  \begin{center}
  { \it 
 THEMA, CY Cergy Paris Université, 33 Boulevard du Port, 95011 Cergy-Pontoise Cedex, France\\
 % E-mail: william.kengne@u-cergy.fr\\
 % $^{\text{b}}$   LERSTAD, Université Gaston Berger, Saint-Louis, Sénégal. \\
  E-mail: mamadou-lamine.diop@cyu.fr ; william.kengne@cyu.fr  \\
  }
\end{center}

 \pagestyle{myheadings}
 \markboth{Statistical learning for $\psi$-weakly dependent processes}{Diop and Kengne}

~~\\
\textbf{Abstract}:
We consider statistical learning question for $\psi$-weakly dependent processes, that unifies a large class of weak dependence conditions such as mixing, association,$\cdots$
The consistency of the empirical risk minimization algorithm is established. 
We derive the generalization bounds and provide the learning rate, which, on some H{\"o}lder class of hypothesis, is close to the usual $O(n^{-1/2})$ obtained in the {\it i.i.d.} case. Application to time series prediction is carried out with an example of causal models with exogenous covariates. 
 \medskip
 

\section{Introduction}

 Statistical learning has received a considerable attention in the literature in the recent decades. This interest is still increasing nowadays, mainly due to the significant successes of the applications of machine learning algorithms and the theoretical properties of these algorithms, which are now well studied in many cases. 
 For instance, see \cite{vapnik1999nature}, \cite{cucker2002mathematical}, \cite{bousquet2003introduction}, \cite{mohri2018foundations} for some results when the training samples are independent and identically
distributed ({\it i.i.d.}). But, the {\it i.i.d.} assumption fails in many real life applications: market prediction, GDP (gross domestic product) prediction, signal processing, meteorological observations,$\cdots$
There is a vast literature on learning with dependent observations, we refer to \cite{steinwart2009learning}, \cite{zou2009generalization}, \cite{steinwart2009fast}, \cite{vidyasagar2013learning}, \cite{kuznetsov2017generalization} and the references therein for an overview of this issue.   
 
\medskip

We consider the supervised learning setting and let $D_n=\{Z_1=(X_1,Y_1),\cdots, Z_n=(X_n,Y_n) \}$ (the training sample) be a trajectory  of a stationary and ergodic process $\{Z_t=(X_t,Y_t), ~ t \in \Z \}$, which takes values in $ \mathcal{Z} = \mathcal{X} \times \mathcal{Y}$, where $\mathcal{X}$ is the input space and $\mathcal{Y}$ the output space. Denote by $\mathcal{H} = \{ h: \mathcal{X} \rightarrow  \mathcal{Y}\}$ the set of hypothesis (a family of predictors) and consider a loss function $\ell : \mathcal{Y} \times \mathcal{Y} \rightarrow [0,\infty)$.
In this context of learning from dependent data, the generalization error (risk) can be defined in various ways, see for instance \cite{kuznetsov2017generalization}. We deal with the widely-used averaged risk (see for instance, \cite{zou2009generalization}, \cite{steinwart2009fast}, \cite{alquier2013prediction},  \cite{hang2014fast}, \cite{blanchard2019concentration}), given for any hypothesis $h \in \mathcal{H}$ by 
\[  R(h) = \E\big[\ell \big(h(X_0),Y_0 \big) \big] .\]
The goal is to construct a learner $h \in \mathcal{H}$ such that, for any $t\in \Z$, $h(X_t)$ is average "close" to $Y_t$; that is, a learner which achieves the small averaged risk. The empirical risk (with respect to the training sample) of a hypothesis $h$ is given by
\[ \widehat{R}_n(h) = \frac{1}{n} \sum_{i=1}^n \ell\big(h(X_i),Y_i \big)   .\] 
In the sequel, we set $\ell(h,z) = \ell\big(h(x),y\big)$ for all $z=(x,y) \in  \mathcal{X} \times\mathcal{Y}$ and $h: \mathcal{X} \rightarrow \mathcal{Y}$.
 The setting considered here covers many commonly used situations: regression estimation, classification (pattern recognition when $\mathcal{Y}$ is finite), autoregressive models prediction (we can take $X_t=(Y_{t-1},\cdots,Y_{t-k})$ and $\mathcal{X}= \mathcal{Y}^k$ for some $k\in \N$), autoregressive models with exogenous covariates.

\medskip

Consider a target (with respect to $\mathcal{H}$) function $h_{\mathcal{H}}$ (assumed to exist), given by,
\[ h_{\mathcal{H}} = \argmin_{h \in \mathcal{H}} R(h);   \] 
and the empirical target 
\begin{equation}\label{def_hat_h}
 \widehat{h}_{n} = \argmin_{h \in \mathcal{H}} \widehat{R}_n(h).
\end{equation}
We focus on the empirical risk minimization (ERM) principle and aim to study the relevance of estimation of $h_{\mathcal{H}}$ by $\widehat{h}_{n}$.
The capacity of $\widehat{h}_{n}$ to approximate $h_{\mathcal{H}}$ is now as the generalization capability of the ERM algorithm.
This generalization capability is accessed by studying how $R(\widehat{h}_{n})$ is close to $R(h_{\mathcal{H}})$. The deviation between $R(\widehat{h}_{n})$ and $R(h_{\mathcal{H}})$ is the generalization error the algorithm.
When $R(\widehat{h}_{n}) - R(h_{\mathcal{H}}) = o_P(1)$, the ERM algorithm is said to be consistent within the hypothesis class $\mathcal{H}$. 

\medskip

The study of a learning algorithm includes the calibration of a bound of the generalization error for any fixed $n$ (non asymptotic property) and the investigation of consistency (asymptotic property).
There exist several important contributions in the literature devoted to statistical learning for dependent observations, with various types of dependence structure.
See among others papers, \cite{mohri2007stability}, \cite{mohri2008rademacher}, \cite{zou2009generalization}, \cite{steinwart2009fast}, \cite{zou2011learning}, \cite{hang2017bernstein}, \cite{kuznetsov2017generalization} for some developments under mixing conditions and \cite{zou2009learning}, \cite{zou2012generalization}, \cite{zou2014generalization}, \cite{xu2014generalization} for some results for  Markov chains.   
\cite{alquier2013prediction} considered a prediction of time series under $\theta$-weakly dependent condition. They established convergence rates using the PAC-Bayesian approach. 
See also \cite{blanchard2019concentration}, \cite{kuznetsov2016time}, \cite{kuznetsov2018theory} for some Bernstein-type inequality for $\tau$-mixing process and some advances on time series forecasting using  statistical learning paradigm. 
However, most of the above works are developed within a mixing condition or for time series prediction or do not consider a general setting that includes pattern recognition, regression estimation, time series prediction,$\cdots$ 

\medskip

In this new contribution, we consider a general learning framework where the observations $D_n=\{Z_1=(X_1,Y_1),\cdots, Z_n=(X_n,Y_n) \}$ is a trajectory of a $\psi$-weakly dependent process $\{Z_t=(X_t,Y_t), ~ t \in \Z \}$ with values in a Banach space $ \mathcal{Z} = \mathcal{X} \times \mathcal{Y}$. The following issues are addressed.

 \begin{enumerate}
 \item[(i)] \textbf{Consistency of the ERM algorithm}. We establish the consistency of the ERM algorithm within any space $\mathcal{H}$ of Lipschitz predictors. In comparison with the existing works, let us stress that, the $\psi$-weakly dependent structure considered here is a more general concept and it is well known that many weak dependent processes do not fulfill the mixing conditions, see for instance \cite{dedecker2007weak}.
\item[(ii)] \textbf{Generalization bounds and convergence rates}. When $\mathcal{X} \subset \R^d$ (with $d \in \N$), $\mathcal{Y} \subset \R$ and $\mathcal H$ is a subset of a H{\"o}lder space $\mathcal{C}^s$ for some $s>0$, generalization bounds are derived and the learning rate is provided. This rate is close to the usual $O(n^{-1/2})$ when $s\gg d$.
\item[(iii)] \textbf{Application to time series prediction}. Application to the prediction of affine causal models with exogenous covariates is carried out. We show that, these models fulfill the conditions that ensure the consistency of the ERM algorithm and enjoy the generalization bounds established. 
 \end{enumerate}

\medskip

The rest of the paper is structured as follows. In Section 2, we set some notations and assumptions. Section 3 provides the main results on consistency, generalization bounds and convergence rates. Application to the prediction of affine causal models with exogenous covariates is carried out in Section 4, whereas Section 5 is devoted to the proofs of the main results.

\section{Notations and assumptions}

In the sequel, we assume that $\mathcal{X}$ and $\mathcal{Y}$ are subsets of separable Banach spaces equipped with norms $\|\cdot\|_{\mathcal{X}}$ and $\|\cdot\|_{\mathcal{Y}}$ respectively and consider the covering number as the complexity measure of $\mathcal{H}$.
The complexity of the hypothesis set $\mathcal{H}$ play a key role in such study.
 Recall that, for any $\epsilon>0$, the $\epsilon$-covering number $\mathcal{N}(\mathcal{H},\epsilon)$ of $\mathcal{H}$ with the $\|\cdot\|_\infty$ norm is the minimal number of balls of radius $\epsilon$ needed to cover  $\mathcal{H}$; that is,
\[ \mathcal{N}(\mathcal{H},\epsilon)= \inf\Big\{ m \geq 1 ~ : \exists h_1, \cdots, h_m \in \mathcal{H} \text{ such that } \mathcal{H} \subset \bigcup_{i=1}^m B(h_i,\epsilon)    \Big\} \]
where 
$B(h_i,\epsilon) = \big\{ h : \mathcal{X} \rightarrow \mathcal{Y}; \| h - h_i\|_\infty = \sup_{x \in \mathcal{X}} \| h(x) - h_i(x)\|_{\mathcal{Y}} \leq \epsilon \big\}$.
For two Banach spaces $E_1$ and $E_2$  equipped with a norm $\|\cdot \|_{E_i}$ $i=1,2$, and any function  $h:E_1 \rightarrow E_2$ set,
\[ \| h\|_\infty = \sup_{x \in E_1} \| h(x) \|_{E_2}, ~ 
\lip_\alpha (h) \coloneqq \underset{x_1, x_2 \in E_1, ~ x_1\neq x_2}{\sup} \dfrac{\|h(x_1) - h(x_2)\|_{E_2}}{\| x_1- x_2 \|^\alpha_{E_1}} 
 \text{ for any } \alpha \in [0,1] ,\]
and for any $\mathcal{K}>0$,
$\Lambda_{\alpha,\mathcal{K}}(E_1,E_2)$ (simply $\Lambda_{\alpha,\mathcal{K}}(E_1)$ when $E_2=\R$) denotes the set of functions  $h:E^u_1 \rightarrow E_2$ for some $u \in \N$, such that  $\|h\|_\infty < \infty$ and  $\lip_\alpha(h) \leq \mathcal{K}$.
When $\alpha=1$, we set  $\lip_1 (h)=\lip(h)$.
 %= \{ h:E\rightarrow \R; ~ \|h\|_\infty < \infty, ~  \lip(h) \leq 1\}$.
%
In the whole sequel, $\Lambda_{1}(E_1)$ denotes $\Lambda_{1,1}(E_1,\R)$ and if $E_1,\cdots,E_k$ are separable Banach spaces equipped with norms $\|\cdot\|_{E_1},\cdots,\|\cdot\|_{E_k}$ respectively, then we set $\|x\|=\|x_1\|_{E_1}+\cdots+\|x_k\|_{E_k}$ for any $x=(x_1,\cdots,x_k) \in E_1 \times \cdots \times E_k$.
%
%Also, we consider the norm $\|(x,y)\|_{\mathcal{Z}}=\|x\|_{\mathcal{X}} + \|y\|_{\mathcal{Y}}$ in the space $\mathcal{Z}$, where $(x,y) \in \mathcal{Z}$. 
%
%
We set the following assumptions for the sequel.
 \begin{enumerate} 
     \item [(\textbf{A1}):] There exists $\mathcal{K}_{\mathcal{H}} >0$ such that $\mathcal{H}$ is a subset of  $ \Lambda_{1,\mathcal{K}_{\mathcal{H}}}(\mathcal{X},\mathcal{Y})$ and  
 $ \sup_{h \in \mathcal{H}} \|h\|_\infty < \infty$. 
  \item [(\textbf{A2}):] There exists $\mathcal{K}_{\ell} >0$ such that, the loss function $\ell$ belongs to  $ \Lambda_{1,\mathcal{K}_{\ell}}(\mathcal{Y}^2)$ and $M = \sup_{h \in \mathcal{H}} \sup_{z \in \mathcal{Z}} |\ell(h,z)| < \infty $.  
  \end{enumerate}   

\medskip

 Under (\textbf{A1}) and (\textbf{A2}), we have 
 \begin{equation} \label{def_L}
  L := \sup_{h_1, h_2 \in \mathcal{H}, h_1 \neq h_2 } \sup_{ z\in \mathcal{Z}  } \frac{ |\ell(h_1,z) - \ell(h_2,z) | }{ \|h_1 - h_2 \|_\infty } < \infty.
 \end{equation}
Under the pre-compact condition in (\textbf{A1}), for any $\epsilon >0$, the $\epsilon$-covering number $\mathcal{N}(\mathcal{H},\epsilon)$ is finite.
If $\mathcal{Y}$ is bounded, then examples of loss functions 
\begin{equation} \label{exampl_loss}
\ell(y, y') = \| y - y' \|_\mathcal{Y}, ~ ~ \ell(y, y') = \| y - y' \|^2_\mathcal{Y},
\end{equation}
fulfill the conditions in assumption (\textbf{A2}) with $\mathcal{K}_{\ell}=1$ and $\mathcal{K}_{\ell}=2 \sup_{y \in \mathcal{Y}} \|y\|_\mathcal{Y}$ respectively. In both cases, one can easily see that $\sup_{h \in \mathcal{H}} \sup_{z \in \mathcal{Z}} |\ell(h,z)| < \infty$.  

\medskip
\noindent
Let us define the weakly dependent process, see \cite{doukhan1999new} and \cite{dedecker2007weak}.
Let $E$ be a separable Banach space.  
%For a set $E$, %$$h:E \rightarrow \R$, 
%we denote by $\Lambda(E)$  the set of functions $h:E \rightarrow \R$ such that $\lip(h)<\infty$ and 
%$\Lambda^{(1)} = \left\{h \in \Lambda,~~ \|h\|_{\infty}< \infty,~\lip(h) < 1 \right\}$, where
%\[
%\lip(h) =\sup \left\{ \frac{|h(x) - h(y)|}{\|x-y\|},~x, y \in E, ~ x\neq y\right\}.
%\]
%
\begin{Def}\label{def_weak_dep}
An $E$-valued process $(Z_t)_{t \in \Z}$ is said to be $(\Lambda_1(E),\psi,\epsilon)$-weakly dependent if there exists a function 
$\psi: [0,\infty)^2 \times \N^2 \rightarrow [0,\infty)$ and a sequence $\epsilon=(\epsilon(r))_{r \in \N}$ decreasing
to zero at infinity such that, 
for any $g_1,g_2 \in \Lambda_1(E)$ with $g_1: E^{u} \rightarrow \R$, $g_2:E^{v}  \rightarrow \R$ ($u, v \in \N$)
 and for any $u$-tuple $(s_1,\cdots,s_u)$ and any $v$-tuple $(t_1,\cdots,t_v)$ with $s_1 \leq \cdots \leq s_u \leq s_u +r \leq t_1 \leq \cdots \leq t_v$,
 the following inequality is fulfilled:
 \begin{equation*}%\label{eq_weak_dep}
 \left|\cov \left(g_1(Z_{s_1},\cdots,Z_{s_u}), g_2(Z_{t_1},\cdots,Z_{t_v})\right) \right| \leq \psi \left(\lip(g_1),\lip(g_2),u,v \right)\epsilon(r).
 \end{equation*}
\end{Def}
 For example, we have the following choices of $\psi$ (see also \cite{dedecker2007weak}).
 \begin{itemize}
 \item $\psi \left(\lip(g_1),\lip(g_2),u,v \right)= v \lip(g_2)$: the $\theta$-weak dependence, then denote $\epsilon(r) = \theta(r)$;
 \item $\psi \left(\lip(g_1),\lip(g_2),u,v \right)= u \lip(g_1) + v \lip(g_2)$: the $\eta$-weak dependence, then denote $\epsilon(r) = \eta(r)$;
 \item $\psi \left(\lip(g_1),\lip(g_2),u,v \right)= u v \lip(g_1) \cdot \lip(g_2)$: the $\kappa$-weak dependence, then denote $\epsilon(r) = \kappa(r)$;
  \item $\psi \left(\lip(g_1),\lip(g_2),u,v \right)= u \lip(g_1) + v \lip(g_2) + u v \lip(g_1) \cdot \lip(g_2)$: the $\lambda$-weak dependence, then denote $\epsilon(r) = \lambda(r)$.   
 \end{itemize}
 One can easily see that $\eta(r) \leq \theta(r)$ for all $r \geq 0$. 
 In the sequel, for each of the four choices of $\psi$ above, we set respectively,
 $\Psi(u,v)=2v$, $\Psi(u,v)=u+v$, $\Psi(u,v)=uv$ and $\Psi(u,v)=(u+v + uv)/2$. 
 
\medskip

Now, we set the weak-dependence assumption.
\begin{enumerate}
 \item [(\textbf{A3}):] Let $\psi: [0,\infty)^2 \times \N^2 \rightarrow [0,\infty)$ be one of the four choices above.
 The process $\{Z_t=(X_t,Y_t), ~ t \in \Z \}$ is stationary ergodic and $(\Lambda_1(\mathcal{Z}),\psi,\epsilon)$-weakly dependent such that, there exist $L_1, L_2, \mu \geq 0$ satisfying 
\begin{equation}\label{cond_under_WD}
\sum_{j\geq 0} (j+1)^k \epsilon_j \leq L_1 L_2^k (k!)^\mu  ~ \text{ for all}  ~  k \geq 0.
\end{equation}   
\end{enumerate}

\medskip

\section{Main results}\label{Sec_results}
\subsection{Generalization bound and consistency}
The following proposition provides an inequality of the deviation of the empirical risk around the risk, for any fixed predictor.
\begin{prop}\label{prop1}
Assume that the conditions  (\textbf{A1})-(\textbf{A3}) hold. Let  $h \in \mathcal{H}$. 
For all $\varepsilon >0$, $n \in \N$,  we have
  
  \begin{equation*}%\label{res_th1}
P\Big \{R(h) - \widehat{R}_n(h) > \varepsilon \Big \} 
\leq
 \exp\Big(- \dfrac{n^2 \varepsilon^2/2}{A_n + B_n^{1/(\mu+2)} (n \varepsilon)^{(2\mu+3)/(\mu+2)}} \Big) ,
\end{equation*}
for any real numbers $A_n$ and $B_n$ satisfying:
 \[
 A_n \geq \E \Big[\Big(  \sum_{i=1}^n \Big(\ell\big(h(X_i),Y_i \big) - \E\big[\ell \big(h(X_0),Y_0 \big) \big]\Big)  \Big)^2 \Big] ~~and~~ B_n=2  M L_2 \max\Big(\dfrac{ 2^{4+\mu} n  M^2 L_1}{A_n}, 1 \Big).
 \]
\end{prop}

\medskip
 
 \noindent The next theorem provides a uniform concentration inequality between risk and its empirical version.
 
 \medskip
 
\begin{thm}\label{th1}
Assume that the conditions of Proposition \ref{prop1} hold. 
 For all $\varepsilon >0$, $n \in \N$, we have
 \begin{equation}\label{prob_R_R_n}
P \Big\{ \sup_{h \in \mathcal H} \Big[R(h) - \widehat{R}_n(h)\Big] >  \varepsilon \Big\}
 \leq 
 \mathcal{N}\Big(\mathcal{H},\frac{\varepsilon }{4L}\Big)\exp\Big(- \dfrac{n^2 \varepsilon^2/8}{A_n + B_n^{1/(\mu+2)} (n \varepsilon /2)^{(2\mu+3)/(\mu+2)}} \Big),
\end{equation} 
where $A_n$, $B_n$ are defined in Proposition  \ref{prop1} and $L$ is defined in (\ref{def_L}).
\end{thm}

\medskip

\medskip

\begin{rmrk}\label{rmk1}
As in Remark 2 in \cite{doukhan2007probability}, we can use (\textbf{A3}) to show that 
\[
 \E \Big[\Big(  \sum_{i=1}^n \Big(\ell\big(h(X_i),Y_i \big) - \E\big[\ell \big(h(X_0),Y_0 \big) \big]\Big)  \Big)^2 \Big]
 \leq
 2nM^2 \Psi(1,1)L_1.
\]
Let  $\varepsilon >0$.
Thus, taking $A_n= 2nM^2 \Psi(1,1)L_1$ (which gives  $B_n=2  M L_2 \max(2^{3+\mu}/\Psi(1,1), 1 )$), the following inequality can be immediately deduced from  Theorem \ref{th1}: 
	\begin{equation}\label{rem_R_hat_R}
P \big\{ \sup_{h \in \mathcal H} \big[R(h) - \widehat{R}_n(h)\big] >  \varepsilon \big\}
 \leq 
 \mathcal{N}\Big(\mathcal{H},\frac{\varepsilon}{4L}\Big)\exp\Big(- \dfrac{n^2 \varepsilon^2/4 }{C_1n + 2C_2^{1/(\mu+2)} (n \varepsilon/2)^{(2\mu+3)/(\mu+2)}} \Big),	
 \end{equation}
where $C_1= 4M^2 \Psi(1,1)L_1~$ and $~C_2= 2  M L_2 \max\Big(\dfrac{ 2^{3+\mu}}{ \Psi(1,1)}, 1 \Big)$.
\end{rmrk}
\medskip

\noindent
Note that, as one can see in the proof of Theorem \ref{th1}, (\ref{prob_R_R_n}) still holds with $\widehat{R}_n(h)- R(h)$; that is,
\begin{equation}\label{prob_R_n_R}
P \Big\{ \sup_{h \in \mathcal H} \Big[ \widehat{R}(h) - R_n(h)\Big] >  \varepsilon \Big\}
 \leq 
 \mathcal{N}\Big(\mathcal{H},\frac{\varepsilon }{4L}\Big)\exp\Big(- \dfrac{n^2 \varepsilon^2/8}{A_n + B_n^{1/(\mu+2)} (n \varepsilon /2)^{(2\mu+3)/(\mu+2)}} \Big).
\end{equation} 
 Therefore, by taking $A_n= 2nM^2 \Psi(1,1)L_1$ as in Remark \ref{rmk1}, we have for all $\varepsilon >0$,
\begin{equation}\label{prob_R_R_n_abs}
P \Big\{ \sup_{h \in \mathcal H} \big|R(h) - \widehat{R}_n(h)\big| >  \varepsilon \Big\}
 \leq 
 2  \mathcal{N}\Big(\mathcal{H},\frac{\varepsilon}{4L}\Big)\exp\Big(- \dfrac{n^2 \varepsilon^2/4 }{C_1n + 2C_2^{1/(\mu+2)} (n \varepsilon /2)^{(2\mu+3)/(\mu+2)}} \Big) \limiten 0 .
\end{equation}
Hence,
\begin{align}
\nonumber 0 \leq R(\widehat{h}_n) -R(h_{\mathcal H}) &= \big(R(\widehat{h}_n) - \widehat{R}_n(\widehat{h}_{n}) \big)+ \big( \widehat{R}_n(\widehat{h}_{n}) - \widehat{R}_n(h_{\mathcal H}) \big) + \big(\widehat{R}_n(h_{\mathcal H}) -R(h_{\mathcal H}) \big) \\
\nonumber &\leq \sup_{h \in \mathcal H}\big|R(h) - \widehat{R}_n(h)\big| + \big( \widehat{R}_n(\widehat{h}_{n}) - \widehat{R}_n(h_{\mathcal H}) \big) + \sup_{h \in \mathcal H} \big|R(h) - \widehat{R}_n(h)\big| \\
%
%\label{ineq_consit_ERM}
%
\nonumber  &\leq 2 \sup_{h \in \mathcal H}\big|R(h) - \widehat{R}_n(h)\big| \limiteproban 0 ,  
\end{align}
where the last inequality above holds since $\widehat{R}_n(\widehat{h}_{n}) - \widehat{R}_n(h_{\mathcal H}) \leq 0$ from (\ref{def_hat_h}).
Thus, whenever the $\epsilon$-covering number of the hypothesis set $\mathcal{H}$ is finite for all $\epsilon >0$, the EMR algorithm (\ref{def_hat_h}) with weakly dependent observations is consistent.

\medskip

\indent Let us stress that, the equation (\ref{prob_R_R_n}) in Theorem \ref{th1} is a non asymptotic result. The following theorem provides an asymptotic result with weaker condition than (\textbf{A3}).
\begin{thm}\label{th2}
  Assume that (\textbf{A1})-(\textbf{A2}) hold and that $(Z_t)_{t \in \Z}$ is $(\Lambda_1(\mathcal{X} \times \mathcal{Y}),\psi,\epsilon)$-weakly dependent with $\epsilon_j =\mathcal{O}(j^{-2})$. 
  For any $\nu \in (0,1]$ and for $n$ large enough, we have for all $\varepsilon >0$,
  \begin{equation}\label{prob_R_R_n_asymp}
P \Big\{ \sup_{h \in \mathcal H} \Big[R(h) - \widehat{R}_n(h)\Big] >  \varepsilon \Big\}
 \leq 
 C_3 \ \mathcal{N}\Big(\mathcal{H},\frac{\varepsilon }{4L}\Big) \exp\Big(\log \log n - \dfrac{n^2 \varepsilon^2/4}{A'_n + B'_n (n \varepsilon/2)^\nu} \Big), 
 \end{equation}
 with $A'_n \geq \E \left[\left(  \sum_{i=1}^n \Big(\ell\big(h(X_i),Y_i \big) - \E\big[\ell \big(h(X_0),Y_0 \big) \big]\Big)  \right)^2 \right]$, 
$B'_n= \dfrac{n^{3/4} \log n}{A'_n} $ and some constant $C_3 >0$.
\end{thm}

\medskip

\noindent Under the assumption of Theorem \ref{th2} with $\epsilon_j =\mathcal{O}(j^{-\gamma})$ for some $\gamma >3$ and from \cite{hwang2013study}, we can find a constant $C>0$ such that for any $h \in \mathcal{H}$,
\begin{equation}\label{bound_var_fast}
\E \Big[\Big(  \sum_{i=1}^n \Big(\ell\big(h(X_i),Y_i \big) - \E\big[\ell \big(h(X_0),Y_0 \big) \big]\Big)  \Big)^2 \Big]
 \leq n C . 
\end{equation}
Therefore, by taking $A'_n= n C$ as in Remark \ref{rmk1},
we can see from Theorem \ref{th2} that $\sup_{h \in \mathcal H}\big|R(h) - \widehat{R}_n(h)\big| = o_P(1)$. This shows the consistency of the EMR algorithm (\ref{def_hat_h}) under a weaker condition than (\textbf{A3}) when the $\epsilon$-covering number of $\mathcal{H}$ is finite for all $\epsilon >0$.

%\medskip

\subsection{Generalization error bounds and convergence rates}

Let us derive the generalization error bounds of the ERM algorithm under the weak dependence conditions.
% In the following proposition,
%
 We deal with the H{\"o}lder space $\mathcal{C}^s$ with $s>0$, and set the assumption,
\begin{enumerate}
 \item [(\textbf{A4}):] $\mathcal{X} \subset \R^d$ (with $d \in \N$) is pre-compact,  $\mathcal{Y} \subset \R$, $\mathcal H$ is a subset of a H{\"o}lder space $\mathcal{C}^s(\overline{\mathcal{X}})$ for some $s>0$, where $\overline{\mathcal{X}}$ denotes the closure  of $\mathcal{X}$. 
 \end{enumerate}
 Let $s>0$. Recall that, the H{\"o}lder space $\mathcal{C}^s(\overline{\mathcal{X}})$ is a set of functions $h: \overline{\mathcal{X}} \rightarrow \R$ such that for any $\beta=(\beta_1,\cdots,\beta_d) \in \N^d$ with $|\beta|=\sum_{i=i}^d \beta_i \leq [s]$, $\| \partial^\beta h \|_\infty < \infty$ and for any  multi-integers $\beta$ with $|\beta| = [s]$, $Lip_{s-[s]}(\partial^\beta h) < \infty$; where $[s]$ denotes the integer part of $s$.
 Equipped with the norm
 \[ \|h \|_{\mathcal{C}^s(\overline{\mathcal{X}})} = \sum_{0\leq |\beta| \leq [s]} \| \partial^\beta h \|_\infty + \sum_{|\beta| = [s]} Lip_{s-[s]}(\partial^\beta h), \]
 $\mathcal{C}^s(\overline{\mathcal{X}})$ is a Banach space (see \cite{Triebel1992}).

 \begin{prop}\label{prop2}
 Assume that (\textbf{A1})-(\textbf{A4}) hold. 
Let $\eta \in (0,1)$. Assume that 
 	\begin{equation}\label{prop_cond_n}
	 n \geq \max \Big\{ \Big(\frac{2C_4}{\sqrt{2M}} \log(1/\eta) \Big)^{\mu+2}, \frac{(2C_0 C_4)^{\mu+2} (4L)^{2(\mu+2)d/s}}{(2M)^{(\mu+2)(2s+2d)/s}} \Big\},
	\end{equation}
 where $C_4 =  4C_1 + 8C_2^{1/(\mu+2)} M^{(2\mu+3)/(\mu+2)}$ (see  Remark \ref{rmk1} for $(C_1, C_2)$) and $C_0$ is given in (\ref{eq_cov_nb}).
\begin{enumerate}
	\item [\rm (i)] With probability at least $1-\eta$, we have
\begin{equation}\label{bound_R_hat_R}
R(\widehat{h}_{n}) - \widehat{R}_n(\widehat{h}_{n}) \leq  \varepsilon_1(n,\eta),
\end{equation}
where 
\[
\varepsilon_1(n,\eta) \leq \max \left\{ \Big[\frac{2\log(1/\eta)}{C_{n,1}} \Big]^{\frac{1}{2}},  \Big[\frac{2C_0[4L]^{\frac{2d}{s}}}{C_{n,1}}\Big]^{\frac{s}{2s+2d}}  \right\}, ~ ~   C_{n,1}=\dfrac{n^2 }{4C_1n + 8C_2^{1/(\mu+2)} (nM)^{(2\mu+3)/(\mu+2)}} .
\]
%
 %$(C_1, C_2)$ is given in Remark \ref{rmk1}.
%
\item [\rm (ii)] With probability at least $1-2\eta$, we have
\begin{equation}\label{bound_R_R}
 R(\widehat{h}_n) -R(h_{\mathcal H}) \leq \varepsilon_1(n,\eta) +\varepsilon'_1(n,\eta),
\end{equation}
where
\begin{flalign*}
& \varepsilon'_1(n,\eta)=\Big[ \frac{\log(1/\eta)}{C_{n,1}}\Big]^{\mu+2}.& &
\end{flalign*}
\end{enumerate}
\end{prop}

In the proposition above, the bound (\ref{bound_R_hat_R}) evaluates the estimation of the risk of $\widehat{h}_n$ from the empirical one, and the bound in (\ref{bound_R_R}) assesses how close this risk to the smallest possible risk on $\mathcal{H}$. The learning rate of the ERM algorithm in this general case is $O(n^{-\frac{s}{(2s+2d)(\mu+2)}})$.
The following proposition provides more faster learning rate.
\begin{prop}\label{prop_fast}
 Assume that (\textbf{A1}), (\textbf{A2}), (\textbf{A4}) hold and that $(Z_t)_{t \in \Z}$ is $(\Lambda_1(\mathcal{X} \times \mathcal{Y}),\psi,\epsilon)$-weakly dependent with $\epsilon_j =\mathcal{O}(j^{-\gamma})$ for some $\gamma >3$. 
Let $\eta \in (0,1)$ and $\nu \in (0,1]$. Assume that 
 	\begin{equation}\label{prop_cond_n_fast}
	n  \geq \max \Big\{ \exp(\eta /C_3), \frac{(C_5 C_3/\eta)^{2}}{4M^4},  \frac{2 C_5 C_0[4L]^{\frac{2d}{s}}}{(2M)^{(2s+2d)/s}} \Big\},
	\end{equation}
 where $C_3$ is given in (\ref{appl_DN07_fast_1}), $C_5 = 4 \big(C + M^\nu /C \big)$ where the constant $C$ is defined in (\ref{bound_var_fast}).
\begin{enumerate}
	\item [\rm (i)] With probability at least $1-\eta$, we have
\begin{equation}\label{bound_R_hat_R_fast}
R(\widehat{h}_{n}) - \widehat{R}_n(\widehat{h}_{n}) \leq  \varepsilon_2(n,\eta,\nu),
\end{equation}
where 
\[
\varepsilon_2(n,\eta,\nu) \leq \max \left\{ \Big[\frac{2\log(C_3 \log n/\eta)}{C_{n,2}} \Big]^{\frac{1}{2}},  \Big[\frac{2C_0[4L]^{\frac{2d}{s}}}{C_{n,2}}\Big]^{\frac{s}{2s+2d}}  \right\}, ~ ~    
C_{n,2}= \dfrac{n^2/4}{ n C + \log n \ n^{\nu-1/4} M^\nu /C}.
\]
\item [\rm (ii)] With probability at least $1-2\eta$, we have
\begin{equation}\label{bound_R_R_fast}
 R(\widehat{h}_n) -R(h_{\mathcal H}) \leq \varepsilon_2(n,\eta,\nu) +\varepsilon_2'(n,\eta,\nu),
\end{equation}
where
\[ \varepsilon'_2(n,\eta,\nu) = \Big[ \frac{\log(C_3 \log n/\eta)}{C'_{n,2}}\Big]^{\frac{1}{2}},~ ~  C'_{n,2}= \dfrac{n^2}{nC +  \log n \ n^{\nu-1/4} (2M)^{\nu}/ (2C)} .\]
\end{enumerate}
\end{prop}
The learning rate in the bounds (\ref{bound_R_hat_R_fast}) and (\ref{bound_R_R_fast}) is $O(n^{-\frac{s}{2s+2d}})$. 
In comparison with existing results based on a H{\"o}lder class, let us note that, this rate is the same to that obtained by \cite{zou2012generalization} and \cite{zou2009learning}  on uniformly and  V-geometrically ergodic Markov chains, and is close to the one from exponentially strongly mixing and geometrically beta-mixing observations, see \cite{zou2009generalization} and \cite{zou2011learning}.   
When $s \gg d$, this learning rate is close to $O(n^{-1/2})$, which is the usual rate obtained in the {\it i.i.d.} case.
%

%%%%%%%%%%%%%%%%%%%%%%%%%%%%%%%%%%%%%%%%%%%%%%%%%%%%%%%%%%%%%%%%%%%%%%%%%
%%%%%%%%%%%%%%%%%%%%%%%%%%%%%%%%%%%%%%%%%%%%%%%%%%%%%%%%%%%%%%%%%%%%%%%%%

\section{Applications to time series prediction}
\subsection{Affine causal models with exogenous covariates}
Let $(\mbox{\large$\chi$}_t)_{t \in \Z}$ be a process of covariates with values in $\R^{d_x}$, $d_{x} \in \N$.  Consider the class of affine causal models with exogenous covariates (see \cite{diop2022inference}) defined by 

 \medskip
 
 \noindent \textbf{Class} $\mathcal {AC}$-$X(\mathcal{M},f):$ A process $\{Y_{t},\,t\in \Z \}$ belongs to $\mathcal {AC}$-$X(\mathcal{M},f)$ if it satisfies:
   \begin{equation}\label{Model_AC_X} 
     Y_t =\mathcal{M}(Y_{t-1}, Y_{t-2}, \ldots; \mbox{\large$\chi$}_{t-1},\mbox{\large$\chi$}_{t-2},\ldots)\xi_t + 
      f(Y_{t-1}, Y_{t-2}, \ldots; \mbox{\large$\chi$}_{t-1},\mbox{\large$\mbox{\large$\chi$}$}_{t-2},\ldots),
   \end{equation}
  where 
   $M,~f :  \R^{\N} \times (\R^{d_x})^{\N} \rightarrow \R$ are two measurable functions,  
   and
   $(\xi_t)_{t \in \Z}$ is a sequence of zero-mean \textit{i.i.d} random variable satisfying $\E(\xi^r_0) < \infty$ for some $r \geq 2$ and $\E(\xi^2_0) =1$. 
   Note that, if $X_t \equiv C$ for some constant $C$ (absence of covariates) then, (\ref{Model_AC_X}) reduces to the classical affine causal models studied among other by \cite{bardet2009asymptotic}, \cite{bardet2012multiple}, \cite{bardet2020consistent}, \cite{kengne2021strongly}. Let us point out that, ARMAX, TARX, GARCH-X, APARCH-X (see \cite{francq2019qml}) models belong to the class $\mathcal {AC}$-$X(\mathcal{M},f)$.
  Other examples of models such as APARCH-X$(\delta,\infty)$, ARX($\infty$)-ARCH($\infty$) introduced by \cite{diop2022inference}, belong to the class $\mathcal {AC}$-$X(\mathcal{M}_\theta,f_\theta)$. 
 We refer also to \cite{diop2022inference}, for the inference for $\mathcal {AC}$-$X(\mathcal{M},f)$, in a semiparametric setting.

  \medskip
 
 To study the stability properties of the model (\ref{Model_AC_X}), we set following Lipschitz-type conditions on the function $f$, $\mathcal{M}$ or $\mathcal{M}^2$, see also \cite{diop2022inference}.   
  In the sequel, $0$ denotes the null vector of any vector space. 
 For $\Psi= f$ or $\mathcal{M}$, consider the assumption,
   \medskip
   
    \noindent 
    \textbf{Assumption} \textbf{A}$ (\Psi)$: $|\Psi(0;0)| <\infty$ and
      there exists two sequences of non-negative real numbers $\big(\alpha_{k,Y}(\Psi) \big)_{k\geq 1} $ and $\big(\alpha_{k,\chi}(\Psi) \big)_{k\geq 1}$ satisfying
     $ \sum\limits_{k=1}^{\infty} \alpha_{k,Y}(\Psi) <\infty$, $ \sum\limits_{k=1}^{\infty} \alpha_{k,\chi}(\Psi) <\infty$;
   such that for any  $(y,x), (y',x') \in \R^{\infty} \times (\R^{d_x})^{\infty}$,
  \[  \big | \Psi_\theta(y;x)- \Psi_\theta(y';x') \big |
  \leq  \sum\limits_{k=1}^{\infty}\alpha_{k,Y}(\Psi) |y_k-y'_k|+\sum\limits_{k=1}^{\infty}\alpha_{k,\chi}(\Psi)\|x_k-x'_k\|, 
  \]
where $\| \cdot\|$ denotes a vector norm in $\R^{q}$, for any $q \in \N$.
\medskip
   
    \noindent
The next assumption is set on the function $H=\mathcal{M}^2$ in the cases of ARCH-X type process.
\medskip
   
  \noindent
 \textbf{Assumption} \textbf{A}$(H)$: Assume that $f=0$, $|\mathcal{M}(0;0)| < \infty$ and there exists two sequences of non-negative real numbers $\big(\alpha_{k,Y}(H) \big)_{k\geq 1} $ and $\big(\alpha_{k,\chi}(H) \big)_{k\geq 1} $ satisfying
     $ \sum\limits_{k=1}^{\infty} \alpha_{k,Y}(H) <\infty $, $ \sum\limits_{k=1}^{\infty} \alpha_{k,\chi}(H) <\infty $;
   such that for any  $(y,x), (y',x') \in \R^{\infty} \times (\R^{d_x})^{\infty}$,
  \[  \big | H(y,x)- H(y',x') \big |
  \leq  \sum\limits_{k=1}^{\infty}\alpha_{k,Y}(H) |y^2_k-{y'}^2_k|+\sum\limits_{k=1}^{\infty}\alpha_{k,\chi}(H) \|x_k-x'_k\| . 
  \]
  
\medskip

 \noindent
 In the sequel, we make the convention that if \textbf{A}$(\mathcal{M})$ holds, then  $\alpha_{k,Y}(H) = \alpha_{k,\chi}(H) =0$ for all $k\geq1$ and if \textbf{A}$(H)$ holds, then $\alpha_{k,Y}(\mathcal{M}) = \alpha_{k,\chi}(\mathcal{M}) =0$ for all $k\geq1$.
 As in \cite{diop2022inference}, we impose an autoregressive-type structure on the covariates:
 \begin{equation}\label{Process_X_AC_X}
 \mbox{\large$\chi$}_t=g(\mbox{\large$\chi$}_{t-1},\mbox{\large$\chi$}_{t-2},\ldots;\eta_t),
\end{equation}
 where $(\eta_t)_{t \in \Z}$ is a sequence of centered random variables with values in $\R^{d_\eta}$ ($d_\eta \in \N$) such as $(\eta_t, \xi_t)_{t \in \Z}$  is {\it  i.i.d}
  and $g$ is a $\R^{d_x}$-valued function such that
\begin{equation} \label{exp_Lip_cov_AC_X}  
   \E\left[\left\|g(0; \eta_0) \right\|^r\right]<\infty
   ~ \text{ and }  ~ 
   \left\|g(x; \eta_0)-g(x'; \eta_0) \right\|_r \leq \sum\limits_{k=1}^{\infty} \alpha_k(g) \left\| x_k-x'_k \right\|
   ~ \text{ for all } x, x'  \in (\R^{d_{x}})^{\infty},
  \end{equation}
   for some $r \geq 1$, a non-negative sequence $(\alpha_k(g))_{k \geq 1}$ satisfying $\sum\limits_{k=1}^{\infty} \alpha_k(g) <1 $; and $\|U\|_r := \big(\E \|U \|^r \big)^{1/r}$ for any random vector $U$.
   
   \medskip
   
   We consider model (\ref{Model_AC_X}) with (\ref{Process_X_AC_X}) and (\ref{exp_Lip_cov_AC_X}), and assume that \textbf{A}$(f)$  and  (\textbf{A}$(\mathcal{M})$  or \textbf{A}$(H)$) hold with 
 \begin{equation} \label{Assum_lip_coefs} 
 \sum\limits_{k=1}^{\infty} \max \left\{\alpha_k(g),\, \alpha_{k,Y}(f) + \|\xi_0\|_r \alpha_{k,Y}(M) + \|\xi_0\|^2_r \alpha_{k,Y}(H) \right\} 
 <1 \text{ for some } r \geq 1; 
 \end{equation}
 %
 %
% or, for ARCH-X type model,
 %
% \begin{equation} \label{Assum_lip_coefs_ARCH_type}
% f\equiv 0, ~ \textbf{A}(H) \text{ holds with }  \|\xi_0\|^2_r  \sum\limits_{k=1}^{\infty} \max \left\{\alpha_k(g),\,\alpha_{k,Y}(H) \right\}
% <1 \text{ for some } r \geq 1.
% \end{equation} 
 %
 %
 Under the conditions (\ref{Assum_lip_coefs}), there exists a $\tau$-weakly dependent stationary, ergodic and non anticipative solution $(Y_t, \mbox{\large$\chi$}_t)_{t\in \Z}$ of (\ref{Model_AC_X}) satisfying $ \|(Y_0, \mbox{\large$\chi$}_0)\|_r <\infty$ (see \cite{diop2022inference}).

\medskip

Let us recall the $\tau$-coefficient (see also \cite{dedecker2004coupling}) for weak dependence.
 Let $(S,\mathcal{A}, \mathbb{P})$ be a probability space, $ \mathcal{V}$ a $\sigma$-subalgebra of $\mathcal{A}$  and $U$ a
random variable with values in a separable Banach space $(E, \| \cdot\|)$,  and satisfying $\|U \|_1 < \infty$. 
Define the $\tau$-coefficient as
\[ \tau(\mathcal{V}, U) = \Big \|\sup_{h\in\Lambda_1(E)}\Big\{\Big|\int h(x) \mathbb{P}_{U|\mathcal{M}}(dx)-\int
h(x) \mathbb{P}_{U}(dx) \Big |\Big \}\Big \|_1.  \]
Now, consider an $E$-valued strictly stationary process $(U_t)_{t \in \Z}$ and set for all $i \in \Z$, $\mathcal{V}_i = \sigma(U_t, ~ t\leq i)$.  
The dependence between the past and the $k$-tuples future of the process $(U_t)_{t\in\Z}$ may be assessed as follows. 
Define
\[  \tau_k(j) =  \max_{1 \leq \ell \leq k}
\frac1\ell \sup\Big\{\tau(\mathcal{V}_i,(U_{j_1},\ldots,U_{j_\ell}))\text{ with } i+j\leq j_1<\cdots <j_\ell \Big\} \text{ and } 
 \tau(j) = \sup_{k>0} \tau_k(j)  .\]
 The process $(U_t)_{t \in \Z}$ is said to be $\tau$-weakly dependent
when $\tau(j)$ tends to 0 as $j \rightarrow \infty$. 
This weak dependence structure implies other dependence notions  such as $\eta$, $\theta$-weak dependence. Indeed, we have for all $j\geq 0$, $\eta(j) \leq \theta(j) \leq \tau (j)$, see Definition \ref{def_weak_dep} and \cite{dedecker2007weak}.
Hence, there exists a solution $(Y_t, \mbox{\large$\chi$}_t)_{t\in \Z}$ of (\ref{Model_AC_X}) which is $\theta$-weakly dependent. Moreover, according to \cite{doukhan2008weakly}, we have as $j \rightarrow \infty$, 
\begin{equation}\label{tau_coef_ACX}
 \theta(j) \leq \tau(j)= O\Big( \inf_{1\leq \iota \leq j} \Big\{ \alpha^{j/\iota} + \sum_{k \geq \iota+1} \alpha_k  \Big \}  \Big),
\end{equation}
where $\alpha_k = \max \big\{\alpha_k(g),\, \alpha_{k,Y}(f) + \|\xi_0\|_r \alpha_{k,Y}(M) + \|\xi_0\|^2_r \alpha_{k,Y}(H) \big\}$  and $\alpha = \sum_{k\geq 1} \alpha_k$.

\subsection{Prediction by learning theory}
 Let $(Y_1, \mbox{\large$\chi$}_1),\cdots,(Y_n, \mbox{\large$\chi$}_n)$ be a trajectory of a stationary and ergodic process $(Y_t, \mbox{\large$\chi$}_t)_{t \in \Z}$ that satisfies (\ref{Model_AC_X}) and (\ref{Process_X_AC_X}).
 The aim is to predict $Y_{n+1}$ from the observations $(Y_n, \mbox{\large$\chi$}_n),\cdots,(Y_{n-p+1}, \mbox{\large$\chi$}_{n-p+1})$, for some fixed $p \in \N$. 
We carry out the learning theory developed in Section \ref{Sec_results} with: $X_t = \big( (Y_{t-1}, \mbox{\large$\chi$}_{t-1}),\cdots,(Y_{t-p}, \mbox{\large$\chi$}_{t-p})  \big)$, $\mathcal{Y} \subset \R$, $\mathcal{X} \subset \big( \R \times \R^{d_x} \big)^p$, $\mathcal{Z} \subset \big( \R \times \R^{d_x} \big)^p \times \R$ and a loss function $\ell: \mathcal{Y} \times \mathcal{Y} \rightarrow [0, \infty)$. 
The set of hypothesis $\mathcal{H}$ is a family of predictors, denoted
\begin{equation}\label{def_H_h_theta}
 \mathcal{H} = \big\{ h_{\theta} : \mathcal{X} \rightarrow \mathcal{Y}, ~ \theta \in \Theta  \big\},
\end{equation}
where $\Theta$ is a compact subset of $\R^d$, for some fixed $d \in \N$. 
A target predictor (assumed to exist) with respect to $\mathcal{H}$ is $h_{\theta_{\mathcal{H}}}$, with
\begin{equation}\label{def_theta_H}
\theta_{\mathcal{H}} = \argmin_{\theta \in \Theta} R(h_\theta);
\end{equation}
and the predictor obtained from the ERM procedure is $h_{\hat{\theta}}$, with
\begin{equation}\label{def_hat_theta}
\hat{\theta} = \argmin_{\theta \in \Theta} \widehat{R}(h_\theta).
\end{equation}
For any $t \in \Z$ and $\theta \in \Theta$, denote,
\begin{equation*}
\hat{Y}_t^{\theta} = h_{\theta}(X_t).
\end{equation*} 
Therefore, the prediction of $Y_{n+1}$ according to the any predictor $h_{\theta}$ is $\hat{Y}_{n+1}^{\theta}$.
\paragraph{Example of linear autoregressive predictors:}  ~ \\
 These predictors are defined for all $\theta \in \Theta$ by,
 \begin{equation}\label{linear_pred}
  h_{\theta}(X_t) = \sum_{j=0}^q \theta_j^T \cdot \phi_j(X_t),
\end{equation}  
 where $q \in \N$ (fixed), $(\phi_j)_{j \geq 0}$ is a family of functions (for example, a wavelet basis, splines,$\cdots$), with  $\phi_j: \big( \R \times \R^{d_x} \big)^p \rightarrow \R^{d_0}$, $\theta_j \in \R^{d_0}$, $\theta=(\theta_0,\theta_1,\cdots,\theta_q) \in \Theta$ and $^T$ denotes the transpose.
 For example, if $d_0=d_x +1$, $q=p$, $\phi_0\equiv (1,0,0,\cdots,0)^T$, $\phi_1(X_t) = (Y_{t-1}, \mbox{\large$\chi$}_{t-1})^T, \cdots, \phi_q(X_t) = (Y_{t-q}, \mbox{\large$\chi$}_{t-q})^T \in \R \times \R^{d_x}$, $\theta_0=(\theta_{0,0}, 0,0, \cdots, 0)$, $\theta_1=(\theta_{1,Y}, \theta_{1,\chi})^T,\cdots, \theta_q=(\theta_{q,Y}, \theta_{q,\chi})^T \in \R \times \R^{d_x}$, then 
  \begin{equation}\label{linear_AR_pred}
  h_{\theta}(X_t) = \theta_{0,0} +  \sum_{j=1}^q \big(\theta_{j,Y} Y_{t-j} + \theta_{j,\chi}^T \cdot \mbox{\large$\chi$}_{t-j} \big).
\end{equation} 

\medskip

Now, consider the learning framework above with a general family of predictor $\mathcal{H}$ at (\ref{def_H_h_theta}), with the assumption that, for all $v_1,\cdots, v_p \in \R \times \R^{d_x}$, the function $\theta \mapsto h_\theta(v_1,\cdots, v_p)$ is continuous on $\Theta$. 
Let us check the assumptions (\textbf{A1})-(\textbf{A4}) for the class $\mathcal {AC}$-$X(\mathcal{M},f)$.
 \begin{enumerate} 
     \item [(i)] Assume that, for any $\theta \in \Theta$, there exists a sequence of non-negative real numbers $(\alpha_k(\theta))_{1\leq k \leq p}$ with $\max_{1\leq k \leq p} \big( \alpha_k(\theta) \big) \leq \mathcal{K}_{\mathcal{H}}$ for some $\mathcal{K}_{\mathcal{H}} >0$, such that for any $v_1,\cdots, v_p, v'_1,\cdots, v'_p \in \R \times \R^{d_x}$,
\begin{equation}\label{Lip_predict}
 |h_\theta(v_1,\cdots ,v_p) - h_\theta(v'_1,\cdots, v'_p) | \leq \sum_{k=1}^p \alpha_k(\theta) \| v_k - v'_k \| .   
\end{equation} 
In addition, if $\mathcal{X}$ is bounded or $\sup_{\theta \in \Theta} \| h_\theta \|_\infty < \infty$, then, (\textbf{A1}) holds. 
Note that, for a solution $\{V_t=(Y_t, \bchi_t), ~ t\in \Z \}$ of (\ref{Model_AC_X})  and (\ref{Process_X_AC_X}), we have 
\[ V_t = F(V_{t-1},\cdots; U_t), \]
where $U_t=(\xi_t, \eta_t)$ and $F(v;U_t) = \Big( \mathcal{M}(y_1,\cdots; \schi_1,\cdots)\xi_t + f(y_1,\cdots; \schi_1,\cdots) , ~ g(\schi_1,\cdots;\eta_t) \Big)$
for all $v=(y_k, \schi_k)_{k \in \N} \in \big(\R \times \R^{d_x} \big)^\infty$. 
Under the condition (\ref{Assum_lip_coefs}) and for a suitable norm $\| \cdot \|_\omega $ on $\R \times \R^{d_x}$, we have $\E\| F(v;U_0) \|^r_\omega < \infty$ for all $v \in \big(\R \times \R^{d_x} \big)^\infty$ and can find a non-negative sequence $(\alpha_k(F))_{k \in \N}$ satisfying $\sum_{k \geq 1} \alpha_k(F) < 1$ such that
\[ \E \| F(v;U_0) - F(v';U_0) \|^r_\omega \leq \sum_{k \geq 1} \alpha_k(F) \| v_k - v'_k \|_\omega \text{ for all } v, v' \in \big(\R \times \R^{d_x} \big)^\infty,\]
see the proof of Proposition 1 in \cite{{diop2022inference}}. 
Therefore, according to the proof of Theorem 3.1 in \cite{doukhan2008weakly} and the proof of Proposition 4.1 in \cite{alquier2012model}, there exists a function $\widetilde{H}$ and a sequence of non-negative real numbers $(\alpha_k(\widetilde{H}))_{k \geq 1}$, satisfying $\sum_{k\geq 1} \alpha_k(\widetilde{H}) < \infty$ such that, 
\[ V_t = \widetilde{H} (U_t, U_{t-1},\cdots), \]
 and 
\[ \| \widetilde{H}(u_1,\cdots) -  \widetilde{H}(u'_1,\cdots) \| \leq \sum_{k \geq 1} \alpha_k(\widetilde{H}) \|u_k - u'_k \|, ~ \text{ for all } (u_k)_{k\geq 1}, (u'_k)_{k\geq 1} \in \big(\R \times \R^{d_\eta} \big)^\infty . \] 
 Hence, if the innovations $(U_t)_{t \in \Z}$ is bounded (i.e, there exists $b>0$ such that $\| U_0 \| \leq b$ a.s.), then one can easily see that process $(Y_t, \bchi_t)_{t \Z}$ is bounded by $\| \widetilde{H}(0,0,\cdots) \| + b \sum_{k\geq 1} \alpha_k(\widetilde{H})$. Thus, $\mathcal{X}$ can be chosen to be bounded and in addition to (\ref{Lip_predict}), (\textbf{A1}) holds.
\item [(ii)] As stressed in the item (i) above, under  (\ref{Assum_lip_coefs}) and  if the innovations $(U_t)_{t \in \Z}$ is bounded, then $(Y_t)_{t \in \Z}$ is bounded and $\mathcal{Y}$ can be chosen to be bounded. Thus, (\textbf{A2}) holds for the examples of loss functions defined in (\ref{exampl_loss}).
\item[(iii)] Under (\ref{Assum_lip_coefs}), the process $(Y_t, \bchi_t)_{t \in \Z}$ is $\theta$-weakly dependent with the coefficients $\theta(j)$ bounded in (\ref{tau_coef_ACX}).
One can easily see that, $(X_t, Y_t)_{t \in \Z}$ is also $\theta$-weakly dependent with the same coefficients $\theta(j)$. Assume that $(Y_0,\bchi_0)$ is bounded, and let us consider:
\begin{itemize}
\item The geometric case. Assume that,
  \[\alpha_{k,Y}(f) + \alpha_{k,Y}(\mathcal{M}) + \alpha_{k,Y}(H) +  \alpha_k(g) = O(a^k) ~ \text{ for some } ~ a \in [0,1) .  \]
From (\ref{tau_coef_ACX}), we obtain (see also \cite{doukhan2008weakly}) 
$\theta(j)\leq\tau(j)=O\Big( \exp\big(-\sqrt{\log(\alpha) \log(a) j} \big) \Big)$, where $\alpha$ is given in (\ref{tau_coef_ACX}).
Thus, (\ref{cond_under_WD}) holds with $\mu=2$ (see Proposition 8 in  \cite{doukhan2007probability}). The condition on the rate of $(\epsilon_j)_{j \in \N}$ in Theorem \ref{th2} and Proposition \ref{prop_fast} also holds.
\item The Riemanian case. Assume that,
  \[\alpha_{k,Y}(f) + \alpha_{k,Y}(\mathcal{M}) + \alpha_{k,Y}(H) +  \alpha_k(g) = O(k^{-\gamma}) ~ \text{ for some } ~ \gamma > 1.  \]
 From (\ref{tau_coef_ACX}), we get (see also \cite{doukhan2008weakly}) 
$\theta(j)\leq\tau(j)=O\Big( \Big(\dfrac{\log j }{j}  \Big)^{\gamma-1} \Big)$. Thus, the condition $\epsilon_j =\mathcal{O}(j^{-2})$ in Theorem \ref{th2} and Proposition \ref{prop_fast} holds with $\gamma >3$.
\end{itemize}
\item [(iv)]  As above, under  (\ref{Assum_lip_coefs}) and  if the innovations $(U_t)_{t \in \Z}$ is bounded, then $(Y_t, \bchi_t)_{t \in \Z}$ is bounded and $\mathcal{X}$ can be chosen to be a bounded subset of $(\R \times \R^{d_x})^p$. Therefore, the example of the set $\mathcal{H}$ of linear autoregressive predictors $h_\theta (\cdot)$, $\theta \in \Theta$ defined in (\ref{linear_AR_pred}) with $q=p$ is pre-compact and is included in $\mathcal{C}^s(\overline{\mathcal{X}})$ for all $s \in \N$. Thus, (\textbf{A4}) holds for such linear predictors.     
 \end{enumerate}

 \subsection{Some simulation results}
 We present some simulation experiments conducted to assess the generalization bound of the ERM algorithm. 
  To this end, we consider the class $\mathcal{H}$ of linear autoregressive predictors defined in (\ref{linear_AR_pred}) with $q=1$ (hence, $\Theta \subset \R^3$). We evaluate the difference  $R(\widehat{h}_n)  - R(h_{\mathcal H})$ in linear and nonlinear dynamic models.

\medskip

\noindent
Now, consider the ARX(2) and Threshold-ARX(1) (TARX) models given by:
\begin{align}
	\label{ARX}  Y_t&= 0.25 Y_{t-1} -0.4 Y_{t-2}  + 0.8 \bchi_t +\xi_t;\\
	\label{TARX} Y_t&= 0.2 \max \left(Y_{t-1},0\right) -0.6 \min \left(Y_{t-1},0\right) +1.5 \bchi_t +\xi_t;
\end{align}
where $(\xi_t)_{t \in Z}$ represents  innovations generated from a standardized uniform distribution $\mathcal{U}[-2,2]$ and $\bchi_t$ follows an 1-dimensional AR(1) process with nonzero mean and generated with a standardized $\mathcal{U}[-2,2]$ innovations.
%\[
%{\large \chi}_t = \omega_0  +  \omega_1{\large \chi}_{t-1} + \eta_t.
%\]

\medskip

\begin{figure}[h!]
\begin{center}
\includegraphics[height=7.5cm, width=16.5cm]{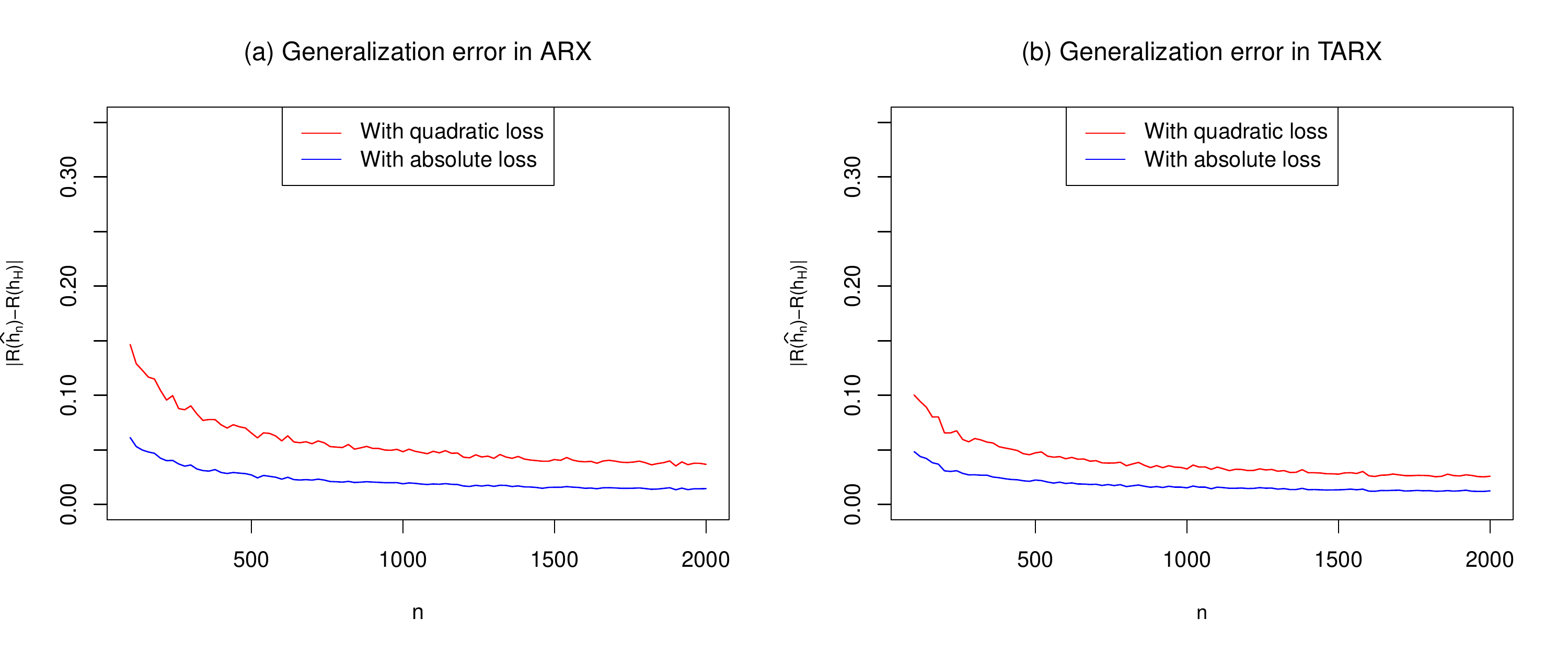}
\end{center}
\vspace{-.7cm}
\caption{\it Plot of $| R(\widehat{h}_n) -R(h_{\mathcal H})|$ obtained with the class of the linear autoregressive predictors (\ref{linear_AR_pred}) 
in the ARX and TARX models.}
\label{Graphe_Sim}
\end{figure}

\noindent
First, for each of the models (\ref{ARX}) and (\ref{TARX}), a sample $(Y_0', \bchi_0'), (Y_1', \bchi_1'), (Y_2', \bchi_2'),\ldots, (Y_m', \bchi_m')$ with $m=10^4$ is generated and $h_{\mathcal H}$ is estimated by $h_{\widetilde{\theta}}$ where
\[ \widetilde{\theta} = \underset{\theta \in \Theta}{\argmin}\widetilde{R}_1(\theta), \text{ with } \widetilde{R}_1(\theta) = \dfrac{1}{m} \sum_{i=1}^m \ell\big( h_\theta(X_i'), Y_i' \big) \text{ and } X_i'=(Y_{i-1}', \bchi_{i-1}').\]
Therefore, $R(h_{\mathcal H})$ is estimated by $\widetilde{R}_1(\widetilde{\theta})$. 
Secondly, for $n=100,120,140,\ldots,2000$, a trajectory $((Y_1, \bchi_1), (Y_2, \bchi_2), \\ \ldots, (Y_n, \bchi_n))$ is generated to compute $\widehat{h}_n$, which is equal to $h_{\widehat{\theta}}$, where 
\[ \widehat{\theta} = \underset{\theta \in \Theta}{\argmin}\widehat{R}_n(h_\theta) .\] 
Hence, $R(\widehat{h}_n)$ is estimated by  $\widetilde{R}_2(\widehat{\theta})$, with 
\[ \widetilde{R}_2(\theta) = \dfrac{1}{n} \sum_{i=1}^n \ell\big( h_\theta(X_i''), Y_i \big) \text{ for all } \theta \in \Theta \text{ and } X_i''=(Y_{i-1}, \bchi_{i-1}), \] 
where $((Y_0, \bchi_0), (Y_1, \bchi_1), (Y_2, \bchi_2),\ldots, (Y_n, \bchi_n))$ is generated from the model and is independent of the sample used for the estimation of $h_{\widetilde{\theta}}$. 
Each simulation for $n$ fixed is repeated 500 times and the Monte Carlo estimate of $R(\widehat{h}_n)$ is computed. 
Figure \ref{Graphe_Sim} displays the curve of points $(n,| R(\widehat{h}_n) -R(h_{\mathcal H}) |)$ obtained with the absolute and squared losses. 
One can see that, for both these losses, the value of $R(\widehat{h}_n) -R(h_{\mathcal H})$ approaching zero when the sample size $n$ increases.
%
%These results are consistent with of  and are overall satisfactory. 
%
These findings are in accordance with the theoretical results of the consistency of the ERM algorithm within the class $\mathcal{H}$ and the generalization bound (\ref{bound_R_R}) in Proposition \ref{prop2} which decreases to zero with high probability. 

%

%%%%%%%%%%%%%%%%%%%%%%%%%%%%%%%%%%%%%%%%%%%%%%%%%%%%%%%%%%%%%%%%%%%%%%%%%%
%%%%%%%%%%%%%%%%%%%%%%%%%%%%%%%%%%%%%%%%%%%%%%%%%%%%%%%%%%%%%%%%%%%%%%%%%%

\section{Proofs of the main results} 
 \subsection{Proof of Proposition \ref{prop1}}
  Let $h \in \mathcal{H}$.  
 Firstly, we derive an exponential inequality for $\widehat{R}_n(h)$. \\
 For any $n\geq 1$, we have
 \[ R(h) - \widehat{R}_n(h) =  \frac{1}{n} \sum_{i=1}^n \Big(\E\big[\ell \big(h(X_0),Y_0 \big) \big]-\ell\big(h(X_i),Y_i \big) \Big) .\]
Since $h \in \Lambda_{1,\mathcal{K}_{\mathcal{H}}}(\mathcal{X},\mathcal{Y})$, $\ell \in \Lambda_{1,\mathcal{K}_\ell}(\mathcal{Y}^2)$ (from (\textbf{A1}) and (\textbf{A2})), then the function $\ell(h(\cdot),\cdot): \mathcal{Z} \rightarrow [0, \infty)$ is $\mathcal{K}_\ell (\mathcal{K}_{\mathcal{H}} + 1)$-Lipschitz.
 Thus, from the hereditary property of the weak dependence (see \cite{dedecker2007weak}), the process $\Big(\ell\big(h(X_t),Y_t \big)\Big)_{t \in \Z}$ is stationary ergodic and weakly dependent.  
Moreover, under (\textbf{A3}), from Proposition 8 and as in Remark 9 in \cite{doukhan2007probability}, one can see that the conditions of Theorem 1 in \cite{doukhan2007probability} are satisfied with $\nu = 0$, $K=M$. 
Let  $(A_n)_{n \in \N}$  a sequence such that $A_n \geq \E \left[\left(  \sum_{i=1}^n \Big(\ell\big(h(X_i),Y_i \big) - \E\big[\ell \big(h(X_0),Y_0 \big) \big]\Big)  \right)^2 \right]$. 
From Theorem 1 of \cite{doukhan2007probability}, we have for all $\varepsilon >0$, 
\begin{align}\label{appl_DN07}
P\big \{ R(h) - \widehat{R}_n(h) > \varepsilon \big \} &= P\Big\{ \sum_{i=1}^n \Big( \E\big[\ell \big(h(X_0),Y_0 \big) \big] -\ell\big(h(X_i),Y_i \big)\Big) >  n \varepsilon \Big\} \nonumber\\
& \leq \exp\Big(- \dfrac{n^2 \varepsilon^2/2}{A_n + B_n^{1/(\mu+2)} (n \varepsilon)^{(2\mu+3)/(\mu+2)}} \Big) 
\end{align}
with $B_n=2  M L_2 \max\Big(\dfrac{ 2^{4+\mu} n  M^2 L_1}{A_n}, 1 \Big)$. 
 Thus, the proposition follows.  \qed 

\subsection{Proof of Theorem \ref{th1}}  
 We proceed as in \cite{zou2009generalization}.\\
 Let $\varepsilon >0$ and $l=\mathcal{N}(\mathcal{H},\frac{\varepsilon}{L})$.
 Denote by $B_j$ (for $j=1,\cdots ,l$) the balls with center  $h_j$
and radius $\varepsilon/L$ needed to cover $\mathcal{H}$; i.e, 
$\mathcal H \subset  B_1 \cup B_2 \cup \cdots \cup  B_l$. 
 We have
 \begin{equation}\label{eq1_p_th2}
 P \Big\{ \sup_{h \in \mathcal H} \bigg[R(h) - \widehat{R}_n(h)\Big] >  \varepsilon \Big\} \leq \sum^l_{j=1} P \Big\{ \sup_{h \in B_j} \Big[R(h) - \widehat{R}_n(h)\Big]  >  \varepsilon \Big\} .
% P \bigg\{ \sup_{h \in \mathcal H} \bigg[\frac{R(h) - \widehat{R}_n(h)}{\sqrt{R(h)}}\bigg] >  \varepsilon \bigg\} \leq \sum^l_{j=1} P \bigg\{ \sup_{h \in B_j} \bigg[\frac{R(h) - \widehat{R}_n(h)}{\sqrt{R(h)}}\bigg]  >  \varepsilon \bigg\} .
  \end{equation}
 Let us define 
 \[
   g(h)= %(1-\delta)
   R(h)-  \widehat{R}_n(h).
 \] 
 For any $h \in B_j$, we get
\begin{align*}
g(h)-g(h_j) &=  [R(h) - R(h_j)] + [\widehat{R}_n(h_j) - \widehat{R}_n(h)] \\   
&\leq   2L\| h-h_j\|_{\infty} 
\leq 2\varepsilon.
\end{align*}
It comes that
\[ \big \{ g(h) \geq 4\varepsilon \big \}  \subset \big \{ g(h_j) \geq  2\varepsilon \big \}  \text{ for any } h \in B_j;\]
which implies
\begin{equation}\label{eq2_p_th2} 
P \Big\{  \sup_{h \in B_j}g(h) \geq 4\varepsilon \Big\} 
\leq
P\Big\{ g(h_j) \geq  2\varepsilon  \Big\}.
\end{equation}
 From (\ref{eq2_p_th2}), 
we have 
\begin{equation*}%\label{eq3_p_th2} 
P \Big\{  \sup_{h \in B_j} \Big[ R(h)- \widehat{R}_n(h) \Big] \geq 4 \varepsilon \Big\} 
\leq
P\Big\{ R(h_j) - \widehat{R}_n(h_j) \geq  2\varepsilon  \Big\}.
\end{equation*}
\medskip
 Then, taking $\varepsilon=  \varepsilon/4$, (\ref{eq1_p_th2}) implies 
 \begin{equation}\label{proof_th1_sup_h_sum_l}
 P \Big\{ \sup_{h \in \mathcal H} \Big[R(h) - \widehat{R}_n(h)\Big] >  \varepsilon \Big\} \leq 
 \sum^l_{j=1} P\Big\{ R(h_j) - \widehat{R}_n(h_j) \geq \varepsilon /2 \Big\}
  \end{equation}
  %
%with $\tau= \varepsilon /2$. 
Hence, according to Proposition \ref{prop1}, we get
\begin{align*}
P \Big\{ \sup_{h \in \mathcal H} \Big[R(h) - \widehat{R}_n(h)\Big] >  \varepsilon \Big\}
& \leq 
 \sum^l_{j=1} \exp\Big(- \dfrac{n^2 \varepsilon^2/8}{A_n + B_n^{1/(\mu+2)} (n \varepsilon /2)^{(2\mu+3)/(\mu+2)}} \Big) \\
 & \leq 
 l\cdot\exp\Big(- \dfrac{n^2 \varepsilon^2/8}{A_n + B_n^{1/(\mu+2)} (n \varepsilon /2)^{(2\mu+3)/(\mu+2)}} \Big)\\
 & \leq 
 \mathcal{N}\Big(\mathcal{H},\frac{\varepsilon }{4L}\Big)\exp\Big(- \dfrac{n^2 \varepsilon^2/8}{A_n + B_n^{1/(\mu+2)} (n \varepsilon /2)^{(2\mu+3)/(\mu+2)}} \Big);
 %& \leq 
 %\mathcal{N}\Big(\mathcal{H},\frac{\tau}{2L}\Big)\exp\Big(- \dfrac{n^2 \tau^2/2}{A_n + B_n^{1/(\mu+2)} (n \tau)^{(2\mu+3)/(\mu+2)}} \Big);
 \end{align*}
 which establishes the  Theorem \ref{th1}.  \qed

\subsection{Proof of Theorem \ref{th2}} 
 From Theorem 2 in \cite{hwang2014note}, we have for all $h \in \mathcal{H}$, $\varepsilon >0$, $\nu \in (0,1]$ and for sufficiently large $n$, 
\begin{align}
P\big \{ R(h) - \widehat{R}_n(h) > \varepsilon \big \} &= P\Big\{ \sum_{i=1}^n \Big( \E\big[\ell \big(h(X_0),Y_0 \big) \big] -\ell\big(h(X_i),Y_i \big)\Big) >  n \varepsilon \Big\} \nonumber\\
&\leq  P\Big\{ \Big|\sum_{i=1}^n \Big( \E\big[\ell \big(h(X_0),Y_0 \big) \big] -\ell\big(h(X_i),Y_i \big)\Big) \Big| >  n \varepsilon \Big\} \nonumber\\
\label{appl_DN07_fast_1}& \leq C_3 \log n \exp\Big(- \dfrac{n^2 \varepsilon^2}{A'_n + B'_n (n \varepsilon)^\nu} \Big)  \\
\label{appl_DN07_fast_2} & \leq C_3 \exp\Big(\log \log n - \dfrac{n^2 \varepsilon^2}{A'_n + B'_n (n \varepsilon)^\nu} \Big), 
\end{align}
for some constant $C_3>0$, any sequence $(A_n')_{n \in \N}$, satisfying $A'_n \geq \E \left[\left(  \sum_{i=1}^n \Big(\ell\big(h(X_i),Y_i \big) - \E\big[\ell \big(h(X_0),Y_0 \big) \big]\Big)  \right)^2 \right]$ and 
$B'_n= \dfrac{n^{3/4} \log n}{A'_n} $.
Set $l = \mathcal{N}\Big(\mathcal{H},\frac{\varepsilon }{4L}\Big)$. From (\ref{proof_th1_sup_h_sum_l}) and (\ref{appl_DN07_fast_2}), we get for $n$ large enough,
\begin{align*}
P \Big\{ \sup_{h \in \mathcal H} \Big[R(h) - \widehat{R}_n(h)\Big] >  \varepsilon \Big\}
& \leq 
 C_3 \sum^l_{j=1} \exp\Big(\log \log n - \dfrac{n^2 \varepsilon^2/4}{A'_n + B'_n (n \varepsilon)^\nu/2} \Big) \\
& \leq 
 C_3 l \cdot \exp\Big(\log \log n - \dfrac{n^2 \varepsilon^2/4}{A'_n + B'_n (n \varepsilon/2)^\nu} \Big) \\
 & \leq 
 C_3 \ \mathcal{N}\Big(\mathcal{H},\frac{\varepsilon }{4L}\Big) \exp\Big(\log \log n - \dfrac{n^2 \varepsilon^2/4}{A'_n + B'_n (n \varepsilon/2)^\nu} \Big). 
\end{align*}
This completes the proof of the theorem.

\subsection{Proof of Proposition \ref{prop2}} 
% 
%
%\begin{enumerate}
%	\item [(i)] 
	 (i) Let $\varepsilon >0$.
	Since $\overline{\mathcal H}$ is compact (from (\textbf{A4})),  there exists a constant $C_0 > 0$ such that (as in \cite{zhou2003capacity} and \cite{zou2009generalization})
\begin{equation}\label{eq_cov_nb}
 \mathcal{N}\Big(\mathcal{H},\frac{\varepsilon }{4L}\Big) \leq \exp \Big( C_0\Big[\frac{\varepsilon}{4L}\Big]^{\frac{-2d}{s}}\Big).
\end{equation}
 Note that, if $\varepsilon > 2M$, then $P \Big\{ \sup_{h \in \mathcal H} \Big[R(h) - \widehat{R}_n(h)\Big] >  \varepsilon \Big\} = 0$. 
For $\varepsilon \in (0,2M]$, according to (\ref{rem_R_hat_R}), we have 
\begin{align}\label{eq1_proof_prop2}
P \Big\{ \sup_{h \in \mathcal H} \Big[R(h) - \widehat{R}_n(h)\Big] >  \varepsilon \Big\}
%& \leq 
% \mathcal{N}\Big(\mathcal{H},\frac{\varepsilon }{4L}\Big) \exp\Big(- \dfrac{n^2 \varepsilon^2/8}{A_n + B_n^{1/(\mu+2)} (n \varepsilon /2)^{(2\mu+3)/(\mu+2)}} \Big)\\
 & \leq 
 \mathcal{N}\Big(\mathcal{H},\frac{\varepsilon}{4L}\Big)\exp\Big(- \dfrac{n^2 \varepsilon^2/4 }{C_1n + 2C_2^{1/(\mu+2)} (n \varepsilon/2)^{(2\mu+3)/(\mu+2)}} \Big) \nonumber\\
% & \leq 
%\exp \Big( C_0\Big[\frac{\varepsilon}{4L}\Big]^{\frac{-2d}{s}} - \dfrac{n^2 \varepsilon^2/8}{A_n n^{(\mu+1)/(\mu+2)} (\varepsilon /2)^{\frac{2\mu+3}{\mu+2}} + B_n^{1/(\mu+2)} (n \varepsilon /2)^{\frac{2\mu+3}{\mu+2}}} \Big) \nonumber\\
& \leq 
\exp \Big( C_0\Big[\frac{\varepsilon}{4L}\Big]^{\frac{-2d}{s}} - \dfrac{n^2 \varepsilon^2/4 }{C_1n + 2C_2^{1/(\mu+2)} (nM)^{(2\mu+3)/(\mu+2)}}  \Big).
 \end{align}
Let $0 < \eta < 1$.
Consider the following equation with respect to $\varepsilon$:
\[
 \exp \Big( C_0\Big[\frac{\varepsilon}{4L}\Big]^{\frac{-2d}{s}} -  \dfrac{n^2 \varepsilon^2/4 }{C_1n + 2C_2^{1/(\mu+2)} (nM)^{(2\mu+3)/(\mu+2)}}  \Big)=\eta. 
\]
It is equivalent to 
\begin{equation}\label{eq2_proof_prop2}
\varepsilon^{2+\frac{2d}{s}} - \frac{\log(1/\eta)}{C_{n,1}} \varepsilon^{\frac{2d}{s}} - \frac{C_0[4L]^{\frac{2d}{s}}}{C_{n,1}} =0
\end{equation}
with
\[ 
C_{n,1}= %\dfrac{n^{(\mu+3)/(\mu+2)}}{\left(A_n   + B_n^{1/(\mu+2)} n\right) 2^{(\mu+3)/(\mu+2)}} 
\dfrac{n^2 }{4C_1n + 8C_2^{1/(\mu+2)} (nM)^{(2\mu+3)/(\mu+2)}} .
\]
According to Lemma 7 in \cite{cucker2002best}, the equation (\ref{eq2_proof_prop2}) has an unique positive solution $\varepsilon_1(n,\eta)$ satisfying
\begin{equation}\label{varep_n_eta}
\varepsilon_1(n,\eta) \leq \max \left\{ \Big[\frac{2\log(1/\eta)}{C_{n,1}} \Big]^{\frac{1}{2}},  \Big[\frac{2C_0[4L]^{\frac{2d}{s}}}{C_{n,1}}\Big]^{\frac{s}{2s+2d}}  \right\}.
\end{equation}
Note that, since $(2\mu+3)/(\mu+2) >1$, we get
\[C_{n,1}  \geq \dfrac{n^2 }{ n^{(2\mu+3)/(\mu+2)} \big(4C_1 + 8C_2^{1/(\mu+2)} M^{(2\mu+3)/(\mu+2)} \big)}  = \frac{n^{\frac{1}{\mu+2}}}{C_4} ,\]
with $C_4 =  4C_1 + 8C_2^{1/(\mu+2)} M^{(2\mu+3)/(\mu+2)} $.
Hence, according to (\ref{prop_cond_n}) and (\ref{varep_n_eta}) ,
 one can easily see that $\varepsilon_1(n,\eta) \leq 2M$.
From (\ref{eq1_proof_prop2}), we deduce that with probability at least $1 - \eta$,
%, the following inequality follows for all $h \in \mathcal H$: 
\begin{equation*}
\sup_{h \in \mathcal H} \Big[R(h) - \widehat{R}_n(h)\Big] \leq  \varepsilon_1(n,\eta);
\end{equation*}
%
%Since $\widehat{h}_{n} \in \mathcal H$ (the function that minimizes the empirical risk $\widehat{R}_n(h)$), we get 
which implies 
\begin{equation}\label{eq3_proof_prop2}
R(\widehat{h}_{n}) - \widehat{R}_n(\widehat{h}_{n}) \leq  \varepsilon_1(n,\eta).
\end{equation}
Hence, the first part of the proposition holds.\\

%\item [(ii)] 
\noindent
(ii) Let us note that, going along similar lines as in Proposition \ref{prop1}, we can establish that the following inequality also follows for all $h \in \mathcal H$:
%\[
%P\Big \{R(h) - \widehat{R}_n(h) > \varepsilon \Big \} 
%\leq
% \exp\Big(- \dfrac{n^2 \varepsilon^2/2}{A_n + B_n^{1/(\mu+2)} (n \varepsilon)^{(2\mu+3)/(\mu+2)}} \Big),~  \text{ for all } h \in \mathcal H.
%\]
%Going along similar lines, we can establish that the following inequality follows for all $h \in \mathcal H$:
\[
P\Big \{\widehat{R}_n(h) - R(h) > \varepsilon \Big \} 
\leq
 \exp\Big(- \dfrac{n^2 \varepsilon^2/2}{A_n + B_n^{1/(\mu+2)} (n \varepsilon)^{(2\mu+3)/(\mu+2)}} \Big).
\]
Thus, taking $A_n= 2nM^2 \Psi(1,1)L_1$ (as in Remark \ref{rmk1}),  with the function $h_{\mathcal H}$, it comes that
\begin{align}\label{eq4_proof_prop2}
P\Big \{\widehat{R}_n(h_{\mathcal H}) - R(h_{\mathcal H}) > \varepsilon \Big \} 
&\leq
 \exp\Big(- \dfrac{n^2 \varepsilon^2/2}{A_n + B_n^{1/(\mu+2)} (n \varepsilon)^{(2\mu+3)/(\mu+2)}} \Big) \nonumber \\
 &\leq
 \exp\Big(- \dfrac{n^2 \varepsilon^2}{C_1n + 2C_2^{1/(\mu+2)} (2nM)^{(2\mu+3)/(\mu+2)}}  \Big) \nonumber \\ 
 %&\leq
 %\exp \Big(- \dfrac{n^{(\mu+3)/(\mu+2)}\varepsilon^\frac{1}{\mu+2}}{2\left(A_n   + B_n^{1/(\mu+2)} n\right)}\Big)  \nonumber \\ 
 &\leq
 \exp \Big(- C'_n \varepsilon^{2}\Big),
\end{align}
where $
 C'_n= \dfrac{n^2 }{C_1n + 2C_2^{1/(\mu+2)} (2nM)^{(2\mu+3)/(\mu+2)}} $. \\
Now, consider the equation (with respect to $\varepsilon$) 
\[  \exp \Big(- C'_n \varepsilon^{2}\Big) = \eta .\]
A solution of this equation is 
\[ \varepsilon'_1(n,\eta) = \Big[ \frac{\log(1/\eta)}{C'_n}\Big]^{\frac{1}{2}}.\]
Thus, from (\ref{eq4_proof_prop2}), we have
\begin{equation*}
\widehat{R}_n(h_{\mathcal H}) - R(h_{\mathcal H}) \leq \varepsilon'_1(n,\eta)
\end{equation*}
with probability at least $1-\eta$. 
Since $\widehat{R}_n(\widehat{h}_n) \leq \widehat{R}_n(h_{\mathcal H})$ (because $ \widehat{h}_n =\underset{h \in \mathcal H}{\text{argmin}} [\widehat{R}_n(h)]$), we deduce %(with probability at least $1-\eta$)
\begin{equation}\label{eq5_proof_prop2}
\widehat{R}_n(\widehat{h}_n)  - R(h_{\mathcal H}) \leq \varepsilon'_1(n,\eta).
\end{equation}
Combining the inequalities  (\ref{eq3_proof_prop2}) and (\ref{eq5_proof_prop2}), with probability at least $1-2\eta$, we obtain 
\[
R(\widehat{h}_n) -R(h_{\mathcal H}) \leq \varepsilon_1(n,\eta) +\varepsilon'_1(n,\eta).
\]
% where $\varepsilon_1(n,\eta)$ satisfies (\ref{varep_n_eta}). 
 This completes the proof of the proposition.  \qed 
%\end{enumerate}
%

\subsection{Proof of Proposition \ref{prop_fast}}
%\begin{enumerate}
%	\item [(i)]
(i) According to (\ref{bound_var_fast}), we can take $A'_n= n C$; which gives $B'_n= \dfrac{\log n}{n^{1/4} C} $.  
Thus, for $\varepsilon \in (0,2M]$,  by using (\ref{prob_R_R_n_asymp}) and (\ref{eq_cov_nb}) and going as in the  proof of Proposition \ref{prop2}, we get for $n$ large enough,
\begin{align}\label{eq1_proof_prop3}
P \Big\{ \sup_{h \in \mathcal H} \Big[R(h) - \widehat{R}_n(h)\Big] >  \varepsilon \Big\}
& \leq 
C_3 \exp \Big( C_0\Big[\frac{\varepsilon}{4L}\Big]^{\frac{-2d}{s}} +\log \log n - \dfrac{n^2 \varepsilon^2/4}{n C +  \log n \ n^{\nu-1/4}  M^\nu /C} \Big).
 \end{align}
Now, consider the following equation with respect to $\varepsilon$:
\[
 C_3\exp \Big( C_0\Big[\frac{\varepsilon}{4L}\Big]^{\frac{-2d}{s}} + \log \log n - \dfrac{n^2 \varepsilon^2/4}{n C +  \log n \ n^{\nu-1/4} M^\nu /C}  \Big)=\eta, 
\]
which is the equivalent to 
\begin{equation*}\label{eq2_proof_prop2_bis}
\varepsilon^{2+\frac{2d}{s}} -  \frac{\log(C_3 \log n/\eta)}{C_{n,2}} \varepsilon^{\frac{2d}{s}} - \frac{C_0[4L]^{\frac{2d}{s}}}{C_{n,2}} =0
\end{equation*}
with 
\begin{align*}%\label{def_Cn}
 C_{n,2}&= %\dfrac{n^{(\mu+3)/(\mu+2)}}{\left(A_n   + B_n^{1/(\mu+2)} n\right) 2^{(\mu+3)/(\mu+2)}} 
\dfrac{n^2/4}{n C + \log n \ n^{\nu-1/4} M^\nu / C}   .
\end{align*}
This equation has an unique positive solution $\varepsilon_2(n,\eta,\nu)$ satisfying (see also  \cite{cucker2002best})
\[
\varepsilon_2(n,\eta,\nu) \leq \max \left\{ \Big[\frac{2\log(C_3 \log n/\eta)}{C_{n,2}} \Big]^{\frac{1}{2}},  \Big[\frac{2C_0[4L]^{\frac{2d}{s}}}{C_{n,2}}\Big]^{\frac{s}{2s+2d}}  \right\}.
\]
Note that, since $\nu \in (0,1]$, for $n$ large enough, we have
\[  C_{n,2} \geq \dfrac{n^2/4}{C'_1n + C'_2 n M^\nu} \geq \frac{n}{C_5} \text{ with }  C_5 = 4 \big(C +  M^\nu / C \big). 
     \]
Moreover,
\[ C_3 \log n/\eta > 1 \Rightarrow n > \exp(\eta /C_3) .\]
Therefore, for $n > \exp(\eta /C_3) $ and from the inequality $\log x < \sqrt{x} \leq x$ for $x>1$, we get,
\[ \Big[\frac{2\log(C_3 \log n/\eta)}{C_{n,2}} \Big]^{\frac{1}{2}} \leq  \Big[\frac{2 C_3 \log n/\eta}{C_{n,2}} \Big]^{\frac{1}{2}} \leq   \Big[\frac{2 C_3 \sqrt{n}/\eta}{C_{n,2}} \Big]^{\frac{1}{2}}  \leq \Big[\frac{2 C_5 C_3 \sqrt{n}/\eta}{n} \Big]^{\frac{1}{2}} = \frac{(2 C_5 C_3/\eta)^{1/2}}{n^{1/4}} .\]
Hence,
\[ \frac{(2 C_5 C_3/\eta)^{1/2}}{n^{1/4}} < 2M \Rightarrow n > \frac{(C_5 C_3/\eta)^{2}}{4M^4} .\]
In the same way,
\[ \Big[\frac{2C_0[4L]^{\frac{2d}{s}}}{C_{n,2}}\Big]^{\frac{s}{2s+2d}} \leq \Big[\frac{2 C_5 C_0[4L]^{\frac{2d}{s}}}{n}\Big]^{\frac{s}{2s+2d}} \leq 2M \Rightarrow n >  \frac{2 C_5 C_0[4L]^{\frac{2d}{s}}}{(2M)^{(2s+2d)/s}} .\]
Then, for $n$ large enough and
\[ n  \geq \max \Big\{ \exp(\eta /C_3), \frac{(C_5 C_3/\eta)^{2}}{4M^4},  \frac{2 C_5 C_0[4L]^{\frac{2d}{s}}}{(2M)^{(2p+2d)/s}} \Big\},\]
it holds that $\varepsilon_2(n,\eta,\nu) \leq 2M$. \\
Thus, using (\ref{eq1_proof_prop3}), we deduce that with probability at least $1 - \eta$,
%, the following inequality follows for all $h \in \mathcal H$: 
%\begin{equation*}
%\sup_{h \in \mathcal H} \Big[R(h) - \widehat{R}_n(h)\Big] \leq  \varepsilon_1(n,\eta);
%\end{equation*}
%Since $\widehat{h}_{n} \in \mathcal H$ (the function that minimizes the empirical risk $\widehat{R}_n(h)$), we get 
%which implies 
\begin{equation}\label{eq2_proof_prop3}
R(\widehat{h}_{n}) - \widehat{R}_n(\widehat{h}_{n}) \leq  \varepsilon_2(n,\eta,\nu).
\end{equation}
Hence, the  part (i) holds.\\

\noindent
%\item[(ii)]
(ii) Using the same arguments, we can go along similar lines as in (\ref{eq5_proof_prop2})
(see proof of Proposition \ref{prop2}(ii)) to establish that,  with probability at least $1 - \eta$, 
\begin{equation}\label{eq3_proof_prop3}
\widehat{R}_n(\widehat{h}_n)  - R(h_{\mathcal H}) \leq \varepsilon'_2(n,\eta,\nu),
\end{equation}
where
\[ \varepsilon'_2(n,\eta,\nu) = \big[ \frac{\log(C_3 \log n/\eta)}{C'_{n,2}}\big]^{\frac{1}{2}}~~ \text{ with } ~~ 
C'_{n,2}= \dfrac{n^2}{nC +  \log n \ n^{\nu-1/4} (2M)^{\nu}/ (2C)} . 
\]
Combining (\ref{eq2_proof_prop3}) and (\ref{eq3_proof_prop3}), with probability at least $1-2\eta$, we get
\[
R(\widehat{h}_n) -R(h_{\mathcal H}) \leq \varepsilon_2(n,\eta,\nu) +\varepsilon'_2(n,\eta,\nu),
\]
 which completes the proof of the proposition.  \qed 
%\end{enumerate}

\bibliographystyle{acm}

\end{document}